\documentclass[a4paper]{amsart}

\usepackage{times,verbatim}

\usepackage{amssymb}
\usepackage{enumerate}
\usepackage{stmaryrd}
\usepackage{graphicx}

\theoremstyle{plain}
	\newtheorem{thm}{Theorem}[section]

\theoremstyle{definition}

\DeclareMathOperator{\Gr}{Gr}
\DeclareMathOperator{\Spec}{Spec}
\DeclareMathOperator{\Proj}{Proj}

\def\A{{\mathbb A}}

\def\F{{\mathbb F}}
\def\N{{\mathbb N}}
\def\P{{\mathbb P}}

\def\Z{{\mathbb Z}}
\def\cO{{\mathcal O}}
\def\cR{{\mathcal R}}
\def\cX{{\mathcal X}}
\def\cY{{\mathcal Y}}

\def\fp{{\mathfrak p}}
\def\Fun{{\F_1}}
\def\Funsq{{\F_{1^2}}}

\def\inv{{\textup{inv}}}
\def\={\equiv}
\def\n={\equiv\hspace{-10,5pt}/\hspace{3,5pt}}
\def\hom{{\textup{hom}}}
\def\SA{{\vphantom{\A}}^{+\!}{\A}}   
\def\SP{{\vphantom{\P}}^{+}{\P}}

\newcommand{\gen}[1]{\langle #1 \rangle}
\newcommand{\bpquot}[2]{#1\!\sslash\!#2}
\newcommand{\bpgenquot}[2]{#1\!\sslash\!\gen{#2}}

\title{\bf Projective geometry for blueprints}

\author{Javier L\'{o}pez Pe\~na}
	\address{Department of Mathematics\\ 
		University College London\\
		25 Gower Street, London WC1E 6BT, United Kingdom}
	\email{jlp@math.ucl.ac.uk}
	\thanks{J. L\'opez Pe\~na's research was supported by MCIM grant MTM2010-20940-C02-01, research group FQM-266 (Junta de Andaluc\'ia) and Max-Planck Institute for Mathematics in Bonn.}

\author{Oliver Lorscheid}
	\address{Department of Mathematics, University of Wuppertal, Gau\ss str.\ 20, 42097 Wuppertal, Germany.}
	\email{lorscheid@math.uni-wuppertal.de}

\date{}

\begin{document}

\begin{abstract}
	In this note, we generalize the $\Proj$--construction from usual schemes to blue schemes. This yields the definition of projective space and projective varieties over a blueprint. In particular, it is possible to descend closed subvarieties of a projective space to a canonical $\Fun$--model. We discuss this in case of the Grassmannian $\Gr(2,4)$.
\end{abstract}

\ \vspace{-0pt}
\maketitle
% {\ \vspace{-150pt}\\ \flushright\tiny\bf Version 0.21\\ \today\\ }\vspace{118pt}

\section{Introduction}
	\label{intro}

Blueprints are a common generalization of commutative (semi)rings and monoids. The associated geometric objects, blue schemes, are therefore a common generalization of usual scheme theory and $\Fun$--geometry (as considered by Kato \cite{Kato94}, Deitmar \cite{Deitmar05} and Connes-Consani \cite{CC09}). The possibility of forming semiring schemes allows us to talk about idempotent schemes and tropical schemes (cf.\ \cite{Mikhalkin10}). All this is worked out in \cite{blueprints1}.

It is known, though not covered in literature yet, that the $\Proj$-construction from usual algebraic geometry has an analogue in $\Fun$-geometry (after Kato, Deitmar and Connes-Consani). In this note we describe a generalization of this to blueprints. In private communication, Koen Thas announced a treatment of $\Proj$ for monoidal schemes (see \cite{Thas12}).

We follow the notations and conventions of \cite{blueprints2}. Namely, all blueprints that appear in this note are proper and with a zero. We remark that the following constructions can be carried out for the more general notion of a blueprint as considered in \cite{blueprints1}; the reason that we restrict to proper blueprints with a zero is that this allows us to adopt a notation that is common in $\Fun$-geometry.

Namely, we denote by $\A^n_B$ the (blue) affine $n$-space $\Spec \bigl(B[T_1,\dotsc,T_n]\bigr)$ over a blueprint $B$. In case of a ring, this does not equal the usual affine $n$-space since $B[T_1,\dotsc,T_n]$ is not closed under addition. Therefore, we denote the usual affine $n$-space over a ring $B$ by $\SA^n_B=\Spec\bigl(B[T_1,\dotsc,T_n]^+\bigr)$. Similarly, we use a superscript ``$+$'' for the usual projective space $\SP^n_B$ and the usual Grassmannian $\Gr(k,n)^+_B$ over a ring $B$.

%%%%%%%%%%%%%%%%%%%%%%%%%%%%%%%%%%%%%%%%%%%%

\section{Graded blueprints and $\Proj$}
	\label{section: graded blueprints and proj}

Let $B$ be a blueprint and $M$ a subset of $B$. We say that $M$ is \emph{additively closed in $B$} if for all additive relations $b\= \sum a_i$ with $a_i\in M$ also $b$ is an element of $M$. Note that, in particular, $0$ is an element of $M$. A \emph{graded blueprint} is a blueprint $B$ together with additively closed subsets $B_i$ for $i\in\N$ such that $1\in B_0$, such that for all $i,j\in\N$ and $a\in B_i$, $b\in B_j$, the product $ab$ is an element of $B_{i+j}$ and such that for every $b\in B$, there are a unique finite subset $I$ of $\N$ and unique non-zero elements $a_i\in B_i$ for every $i\in I$ such that $b\=\sum a_i$. An element of $\bigcup_{i\geq 0}B_i$ is called \emph{homogeneous}. If $a\in B_i$ is non-zero, then we say, more specifically, that $a$ is \emph{homogeneous of degree $i$}.

We collect some immediate facts for a graded blueprint $B$ as above. The subset $B_0$ is multiplicatively closed, i.e. $B_0$ can be seen as a subblueprint of $B$. The subblueprint $B_0$ equals $B$ if and only if for all $i>0$, $B_i=\{0\}$. In this case we say that $B$ is \emph{trivially graded}. By the uniqueness of the decomposition into homogeneous elements, we have $B_i\cap B_j=\{0\}$ for $i\neq j$. This means that the union $\bigcup_{i\geq 0}B_i$ has the structure of a wedge product $\bigvee_{i\geq 0} B_i$. Since $\bigvee_{i\geq 0} B_i$ is multiplicatively closed, it can be seen as a subblueprint of $B$. We define $B_\hom=\bigvee_{i\geq 0} B_i$ and call the subblueprint $B_\hom$ the \emph{homogeneous part of $B$}. 

Let $S$ be a multiplicative subset of $B$. If $b/s$ is an element of the localization $S^{-1}B$ where $f$ is homogeneous of degree $i$ and $s$ is homogeneous of degree $j$, then we say that $b/s$ is a homogeneous element of degree $i-j$. We define $S^{-1}B_0$ as the subset of homogeneous elements of degree $0$. It is multiplicatively closed, and inherits thus a subblueprint structure from $S^{-1}B$. If $S$ is the complement of a prime ideal $\fp$, then we write $B_{(\fp)}$ for the subblueprint $(B_\fp)_0$ of homogeneous elements of degree $0$ in $B_\fp$.

An ideal $I$ of a graded blueprint $B$ is called \emph{homogeneous} if it is generated by homogeneous elements, i.e.\ if for every $c\in I$, there are homogeneous elements $p_i,q_j\in I$ and elements $a_i,b_j\in B$ and an additive relation $\sum a_ip_i+c\=\sum b_jq_j$ in $B$. 

Let $B$ be a graded blueprint. Then we define $\Proj B$ as the set of all homogeneous prime ideals $\fp$ of $B$ that do not contain $B_\hom^+=\bigvee_{i>0}B_i$. The set $X=\Proj B$ comes together with the topology that is defined by the basis 
\[
 U_h \quad = \quad \{ \ \fp\in X\ | \ h\notin \fp \ \}
\]
where $h$ ranges through $B_\hom$ and with a structure sheaf $\cO_X$ that is the sheafification of the association $U_h\mapsto B[h^{-1}]_0$ where $B[h^{-1}]$ is the localization of $B$ at $S=\{h^i\}_{i\ge q0}$.

Note that if $B$ is a ring, the above definitions yield the usual construction of $\Proj B$ for graded rings. In complete analogy to the case of graded rings, one proves the following theorem.

\begin{thm}
	The space $X=\Proj B$ together with $\cO_X$ is a blue scheme. The stalk at a point $\fp\in \Proj B$ is $\cO_{x,\fp}=B_{(\fp)}$. If $h\in B_\hom^+$, then $U_h\simeq\Spec B[h^{-1}]_0$. The inclusions $B_0\hookrightarrow B[h^{-1}]_0$ yield morphisms $\Spec B[h^{-1}]_0\to \Spec B_0$, which glue to a structural morphism $\Proj B \to \Spec B_0$. \qed
\end{thm}

If $B$ is a graded blueprint, then the associated semiring $B^+$ inherits a grading. Namely, let $B_\hom=\bigvee_{i\geq 0}B_i$ the homogeneous part of $B$. Then we can define $B_i^+$ as the additive closure of $B_i$ in $B^+$, i.e.\ as the set of all $b\in B$ such that there is an additive relation of the form $b\=\sum a_k$ in $B$ with $a_k\in B_i$. Then $\bigvee B_i^+$ defines a grading of $B^+$. Similarly, the grading of $B$ induces a grading on a tensor product $B\otimes_CD$ with respect to blueprint morphisms $C\to B$ and $C\to D$ under the assumption that the image of $C\to B$ is contained in $B_0$. 
Consequently, a grading of $B$ implies a grading of $B_\inv=B\otimes_\Fun\Funsq$ and of the ring $B^+_\Z=B_\inv^+$. Along the same lines, if both $B$ and $D$ are graded and the images of $C\to B$ and $C\to D$ lie in $B_0$ and $C_0$ respectively, then $B\otimes_C D$ inherits a grading obtained from the gradings of $B$ and $D$.

%%%%%%%%%%%%%%%%%%%%%%%%%%%%%%%%%%%%%%%%%%%%%

\section{Projective space}
	\label{section: projective space}

The functor $\Proj$ allows the definition of the projective space $\P^n_B$ over a blueprint $B$. Namely, the free blueprint $C=B[T_0,\dotsc,T_n]$ over $B$ comes together with a natural grading (cf.\ \cite[Section 1.12]{blueprints1} for the definition of free blueprints). Namely, $C_i$ consists of all monomials $bT_0^{e_0}\dotsb T_n^{e_n}$ such that $e_0+\dotsb+e_n=i$ where $b\in B$. Note that $C_0=B$ and $C_\hom=C$. The projective space $\P^n_B$ is defined as $\Proj B[T_0,\dotsc,T_n]$. It comes together with a structure morphism $\P^n_B\to\Spec B$.

In case of $B=\Fun$, the projective space $\P^n_\Fun$ is the monoidal scheme that is known from $\Fun$-geometry (see \cite{Deitmar08}, \cite[Section 3.1.4]{CLR10}) and \cite[Ex.\ 1.6]{blueprints2}). The topological space of $\P^n_\Fun$ is finite. Its points correspond to the homogeneous prime ideals $(S_i)_{i\in I}$ of $\Fun[S_0,\dotsc,S_n]$ where $I$ ranges through all proper subsets of $\{0,\dotsc,n\}$.

In case of a ring $B$, the projective space $\P^n_B$ does not coincide with the usual projective space since the free blueprint $B[S_0,\dotsc,S_n]$ is not a ring, but merely the blueprint of all monomials of the form $bS_0^{e_0}\dotsb S_n^{e_n}$ with $b\in B$. However, the associated scheme $\SP^n_B=(\P^n_B)^+$ coincides with the usual projective space over $B$, which equals $\Proj B[S_0,\dotsc,S_n]^+$.

%%%%%%%%%%%%%%%%%%%%%%%%%%%%%%%%%%%%%%%%%%%%

\section{Closed subschemes}
\label{section: closed subschemes}

Let $\cX$ be a scheme of finite type. By an \emph{$\Fun$-model of $\cX$} we mean a blue scheme $X$ of finite type such that $X_\Z^+$ is isomorphic to $\cX$. Since a finitely generated $\Z$-algebra is, by definition, generated by a finitely generated multiplicative subset as a $\Z$-module, every scheme of finite type has an $\Fun$-model. It is, on the contrary, true that a scheme of finite type possesses a large number of $\Fun$-models.

Given a scheme $\cX$ with an $\Fun$-model $X$, we can associate to every closed subscheme $\cY$ of $\cX$ the following closed subscheme $Y$ of $X$, which is an $\Fun$-model of $\cY$. In case that $X=\Spec B$ is the spectrum of a blueprint $B=\bpquot A\cR$, and thus $\cX\simeq\Spec B_\Z^+$ is an affine scheme, we can define $Y$ as $\Spec C$ for $C=\bpquot A{\cR(Y)}$ where $\cR(Y)$ is the pre-addition that contains $\sum a_i\=\sum b_j$ whenever $\sum a_i=\sum b_j$ holds in the coordinate ring $\Gamma\cY$ of $\cY$. This is a process that we used already in \cite[Section 3]{blueprints2}.

Since localizations commute with additive closures, i.e.\ $(S^{-1}B)^+_\Z = S^{-1}(B^+_\Z)$ where $S$ is a multiplicative subset of $B$, the above process is compatible with the restriction to affine opens $U\subset X$. This means that given $U=\Spec(S^{-1}B)$, which is an $\Fun$-model for $\cX'=U_\Z^+$, then the $\Fun$--model $Y'$ that is associated to the closed subscheme $\cY' = \cX'\times_\cX\cY$ of $\cX'$ by the above process is the spectrum of the blueprint $S^{-1}C$. Consequently, we can associate with every closed subscheme $\cY$ of a scheme $\cX$ with an $\Fun$-model $X$ a closed subscheme $Y$ of $X$, which is an $\Fun$--model of $\cY$; namely, we apply the above process to all affine open subschemes of $\cX$ and glue them together, which is possible since additive closures commute with localizations.

In case of a projective variety, i.e.\ a closed subscheme $\cY$ of a projective space $\SP^n_\Z$, we derive the following description of the associated $\Fun$-model $Y$ in $\P^n_\Fun$ by homogeneous coordinate rings. Let $C$ be the homogeneous coordinate ring of $\cY$, which is a quotient of $\Z[S_0,\dotsc,S_n]^+$ by a homogeneous ideal $I$. Let $\cR$ be the pre-addition on $\Fun[S_0,\dotsc,S_n]$ that consists of all relations $\sum a_i\=\sum b_j$ such that $\sum a_i=\sum b_j$ in $C$. Then $B=\bpquot{\Fun[S_0,\dotsc,S_n]}{\cR}$ inherits a grading from $\Fun[S_0,\dotsc,S_n]$ by defining $B_i$ as the image of $\Fun[S_0,\dotsc,S_n]_i$ in $B$. Note that $B\subset C$ and that the sets $B_i$ equal the intersections $B_i=C_i\cap B$ for $i\geq 0$ where $C_i$ is the homogeneous part of degree $i$ of $C$. Then the $\Fun$-model $Y$ of $\cY$ equals $\Proj B$.

\section{$\Fun$--models for Grassmannians}
	\label{section: grassmannians}

One of the simplest examples of projective varieties that is not a toric variety (and in particular, not a projective space) is the Grassmann variety $\Gr(2,4)$. The problem of finding models over $\Fun$ for Grassmann varieties was originally posed by Soul\`e in \cite{Soule04}, and solved by the authors by obtaining a torification from the Schubert cell decomposition (cf. \cite{LL09, LL09b}).

In this note, we present $\Fun$-models for Grassmannians as projective varieties defined through (homogeneous) blueprints. The proposed construction for the Grassmannians fits within a more general framework for obtaining blueprints and totally positive blueprints from cluster data (cf. the forthcoming preprint \cite{LopezUN}).

Classically, the homogeneous coordinate ring for the Grassmannian $\Gr(k,n)$ is obtained by quotienting out the homogeneous coordinate ring of the projective space $\mathbb{P}^{\binom{n}{k}-1}$ by the homogeneous ideal generated by the Pl\"ucker relations. A similar construction can be carried out using the framework of (graded) blueprints. In what follows, we make that construction explicit for the Grassmannian $\Gr(2,4)$.

Define the blueprint $\mathcal{O}_{\Fun}(\Gr(2,4))= \bpquot{\Fun[x_{12}, x_{13}, x_{14}, x_{23}, x_{24}, x_{34}]}{\cR}$ where the congruence $\cR$ is generated by the Pl\"ucker relation $x_{12}x_{34} + x_{14}x_{23} \equiv x_{13}x_{24}$\ (the signs have been picked to ensure that the totally positive part of the Grassmannian is preserved, cf.\ \cite{LopezUN}). Since $\cR$ is generated by a homogeneous relation, $\mathcal{O}_{\Fun}(\Gr(2,4))$ inherits a grading from the canonical morphism
\[
	\pi: \Fun[x_{12}, x_{13}, x_{14}, x_{23}, x_{24}, x_{34}] \longrightarrow \bpquot{\Fun[x_{12}, x_{13}, x_{14}, x_{23}, x_{24}, x_{34}]}{\cR}.
\] 
Let $\Gr(2,4)_\Fun := \Proj(\mathcal{O}_{\Fun}(\Gr(2,4)))$. The base extension $\Gr(2,4)^{+}_\Z$ is the usual Grassmannian, and $\pi$ defines a closed embedding of $\Gr(2,4)_\Fun$ into $\mathbb{P}^5_\Fun$, which extends to the classical Pl\"ucker embedding $\Gr(2,4)_\Z^+ \hookrightarrow \SP_\Z^5$. 

Homogeneous prime ideals in $\mathcal{O}_\Fun(Gr(2,4))$ are described by their generators as the proper subsets $I \subsetneq \{x_{12}, x_{13}, x_{14}, x_{23}, x_{24}, x_{25}\}$ such that $I$ is either contained in one of the sets $\{x_{12}, x_{34}\}$, $\{x_{14}, x_{23}\}$, $\{x_{13}, x_{24}\}$, or otherwise $I$ has a nonempty intersection with all three of them. In other words, $I$ cannot contain elements in two of the above sets without also containing an element of the third one.

\begin{figure}[h] 
\centerline{\includegraphics{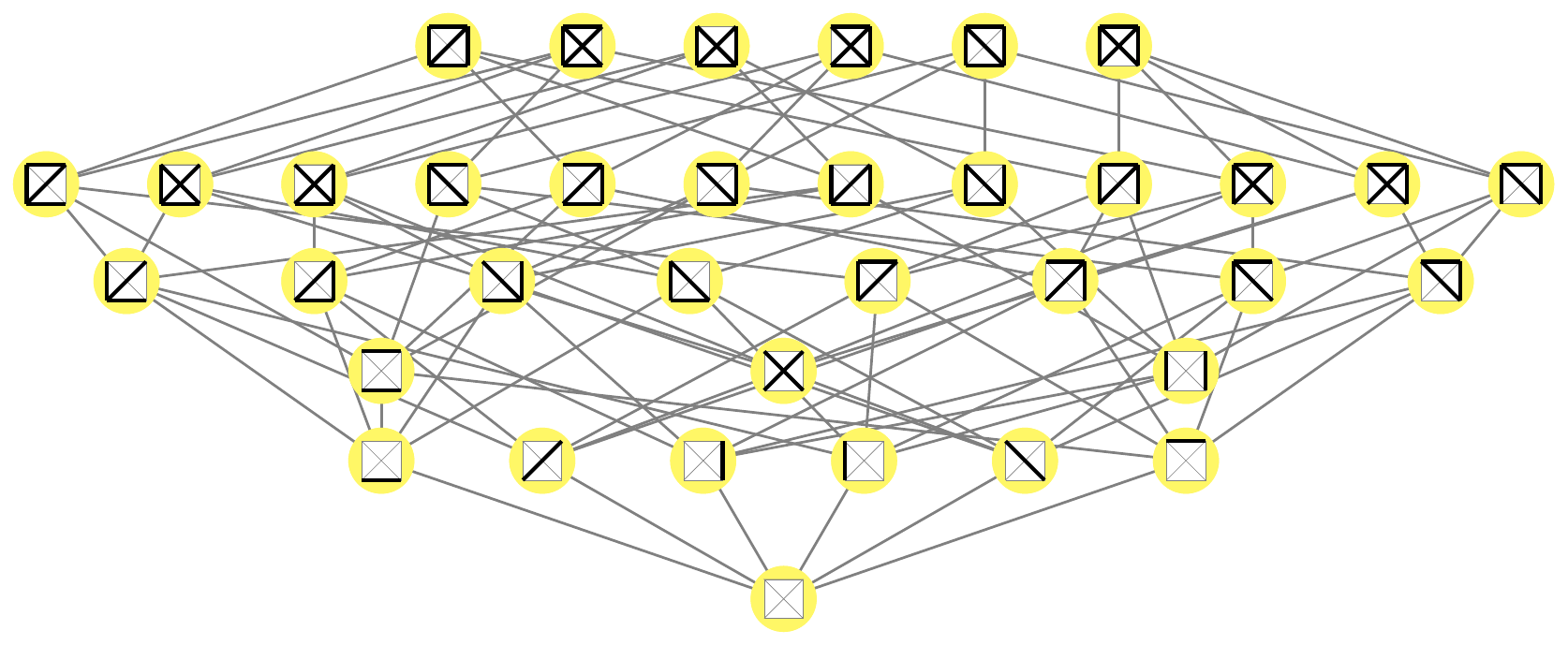}}
\caption{Points of the Grassmannian $\Gr(2,4)_\Fun$. }\label{fig:gr24}
Generator $x_{ij}$ belonging to an ideal is depicted as segment $i$--$j$  in $\vcenter{\hbox{\includegraphics{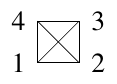}}}$
\end{figure}

The structure of the set of (homogeneous) prime ideals of $\mathcal{O}_\Fun(\Gr(2,4))$ is depicted in Figure \ref{fig:gr24}. It consists of $6 + 12 + 11 + 6 + 1= 36$ prime ideals of ranks $0$, $1$, $2$, $3$ and $4$, respectively (cf.\ \cite[Def.\ 2.3]{blueprints2} for the definition of the rank of a prime ideal), thus resulting in a model essentially different to the one presented in \cite{LL09} by means of torifications, which had $6 + 12 + 11 + 5 + 1 = 35$ points, in correspondence with the coefficients of the counting polynomial $N_{\Gr(2,4)}(q) = 6+12(q-1)+11(q-1)^2+5(q-1)^3+1(q-1)^4$. It is worth noting that despite arising from different constructions, both $\Fun$-models for $\Gr(2,4)$ have $6 = \binom{4}{2}$  closed points, corresponding to the combinatorial interpretation of $\Gr(2,4)_\Fun$ as the set of all subsets with two elements inside a set with four elements. These six points correspond to the $\Fun$-rational Tits points of $\Gr(2,4)_\Fun$, which reflect the naive notion of $\Fun$-rational points of an $\Fun$-scheme (cf.\ \cite[Section 2.2]{blueprints2}).

Like in the classical geometrical setting, the Grassmannian $\Gr(2,4)_{\Fun}$ does admit a covering by six $\Fun$-models of affine $4$-space, which correspond to the open subsets of $\Gr(2,4)_{\Fun}$ where one of $x_{12}$, $x_{34}$, $x_{14}$, $x_{23}$, $x_{13}$ or $x_{24}$ is non-zero. However, these $\Fun$-models of affine $4$-space are not the standard model $\A^4_\Fun=\Spec\bigl(\Fun[a,b,c,d]\bigr)$, but the ``$2\times 2$-matrices'' $M_{2,\Fun}=\Spec\bigl(\bpgenquot{\Fun[a,b,c,d]}{ad\=bc+D}\bigr)$ in case that one of $x_{12}$, $x_{34}$, $x_{14}$ or $x_{23}$ is non-zero, and the ``twisted $2\times 2$-matrices'' $M_{2,\Fun}^\tau=\Spec\bigl(\bpgenquot{\Fun[a,b,c,d]}{ad+bc\=D}\bigr)$ in case that one of $x_{13}$ or $x_{24}$ is non-zero.

\begin{comment}
\newpage
\begin{figure}[h] 
\includegraphics[angle=90,width=0.4\textheight]{Gr24-figure0}
\caption{Points of the Grassmannian $\Gr(2,4)_\Fun$. }\label{fig:gr24}
Generator $x_{ij}$ belonging to an ideal is depicted as segment $i$--$j$  in
$\vcenter{\hbox{\includegraphics{Gr24-figure1}}}$
\end{figure}
\newpage
\end{comment}
 
\begin{small}

\bibliographystyle{plain}

\begin{thebibliography}{10}

\bibitem{CLR10}
Chenghao Chu, Oliver Lorscheid, and Rekha Santhanam.
\newblock Sheaves and {$K$}-theory for {$\mathbb{F}_1$}-schemes.
\newblock {\em Adv. Math.}, 229(4):2239--2286, 2012.

\bibitem{CC09}
Alain Connes and Caterina Consani.
\newblock Characteristic $1$, entropy and the absolute point.
\newblock Preprint, {\tt arXiv:0911.3537v1}, 2009.

\bibitem{Deitmar05}
Anton Deitmar.
\newblock Schemes over {$\mathbb F\sb 1$}.
\newblock In {\em Number fields and function fields---two parallel worlds},
  volume 239 of {\em Progr. Math.}, pages 87--100. Birkh\"auser Boston, Boston,
  MA, 2005.

\bibitem{Deitmar08}
Anton Deitmar.
\newblock {$\mathbb F\sb 1$}-schemes and toric varieties.
\newblock {\em Beitr\"age Algebra Geom.}, 49(2):517--525, 2008.

\bibitem{Kato94}
Kazuya Kato.
\newblock Toric singularities.
\newblock {\em Amer. J. Math.}, 116(5):1073--1099, 1994.

\bibitem{LopezUN}
Javier L\'opez Pe\~na.
\newblock {$\Fun$}-models for cluster algebras and total positivity.
\newblock In preparation.

\bibitem{LL09b}
Javier L\'opez Pe\~na and Oliver Lorscheid.
\newblock Mapping {$\Fun$}-land. an overview of geometries over the field with
  one element.
\newblock In {\em Noncommutative geometry, arithmetic and related topics},
  pages 241--265. Johns Hopkins University Press, 2011.

\bibitem{LL09}
Javier L{\'o}pez~Pe{\~n}a and Oliver Lorscheid.
\newblock Torified varieties and their geometries over {$\mathbb F_1$}.
\newblock {\em Math. Z.}, 267(3-4):605--643, 2011.

\bibitem{blueprints1}
Oliver Lorscheid.
\newblock The geometry of blueprints. {P}art {I}: algebraic background and
  scheme theory.
\newblock {\em Adv. Math.}, 229(3):1804--–1846, 2012.

\bibitem{blueprints2}
Oliver Lorscheid.
\newblock The geometry of blueprints. {P}art {II}: {T}its-{W}eyl models of
  algebraic groups.
\newblock Preprint, {\tt arXiv:1201.1324}, 2012.

\bibitem{Mikhalkin10}
Grigory Mikhalkin.
\newblock Tropical geometry.
\newblock Unpublished notes, 2010.

\bibitem{Soule04}
Christophe Soul{\'e}.
\newblock Les vari\'et\'es sur le corps \`a un \'el\'ement.
\newblock {\em Mosc. Math. J.}, 4(1):217--244, 312, 2004.

\bibitem{Thas12}
Koen Thas.
\newblock Notes on {$\mathbb{F}_1$}, I. Combinatorics of {$\mathcal{D}_0$}-schemes and {$\mathbb{F}_1$}-geometry.
\newblock In preparation.

\end{thebibliography}

\end{small}

\end{document}